\documentclass[12pt]{amsart}
\usepackage[utf8]{inputenc}
\usepackage[english]{babel}
\usepackage{subcaption}
\usepackage{enumerate}
\usepackage{amsfonts}
\usepackage{comment}
\usepackage{amsmath,amssymb}
\usepackage[makeroom]{cancel}
\usepackage{stmaryrd}
\usepackage{tikz-cd} 
\usepackage{enumerate}
\usepackage{enumitem}
\usepackage{graphicx}
\usepackage{amsthm}
\usepackage{wrapfig}
\usepackage{chngcntr}
\counterwithin{figure}{section}
\newtheorem{thm}{Theorem}
\newtheorem*{thm*}{Theorem}
\newtheorem{prop}[thm]{Proposition}

\newtheorem{lemm}[thm]{Lemma}
\theoremstyle{definition}
\newtheorem{defi}{Definition}
\newtheorem*{defi*}{Definition}

\newtheorem{thdef}[thm]{Theorem-Definition}

\counterwithout{figure}{section}

\renewcommand{\int}[1]{\overset{\circ}{#1}}
\newcommand{\clos}[1]{\overline{#1}}

\numberwithin{equation}{section}
\begin{document}
\title{Markov partitions for non-transitive expansive flows} 
\author[I. Iakovoglou]{Ioannis Iakovoglou}
\address{Université Sorbonne Paris Nord, \\ \vspace{0.1cm} UMR 7539 du
CNRS, 99 avenue J.B. Clément, 
93430 Villetaneuse, 
FRANCE.}
\email{e.mail: ioannis.iakovoglou@u-bourgogne.fr}
\begin{abstract}
    In this note we construct Markov partitions for non-transitive expansive flows in dimension 3.
\end{abstract}
\maketitle
According to a classical result of Ratner (see \cite{Ratner}), every transitive Anosov flow admits a Markov partition. A similar result was shown to be true for transitive pseudo-Anosov flows by Brunella in \cite{Brunella}. Historically, Markov partitions have formed  very useful tools for the proof of several classical results in the theory of Anosov and pseudo-Anosov flows. For example, thanks to Ratner's and Brunella's theorem: 
\begin{itemize}
    \item it was shown that a transitive pseudo-Anosov flow is semi-conjugated to the suspension of a shift of finite type
    \item Fried and Brunella showed respectively that any transitive Anosov flow and pseudo-Anosov flow in dimension 3 admits Birkhoff sections (see \cite{Fried} and \cite{Brunella})
\end{itemize}

More recently, in our preprint \cite{preprint}, we used Markov partitions in order to define a new approach to the problem of classification of transitive Anosov flows on $3$-manifolds. During my recent talks, I've been asked several times whether the classification approach developed in \cite{preprint} applies to non-transitive Anosov or pseudo-Anosov flows in dimension 3. An important impediment to such a generalization of \cite{preprint}, was the existence of Markov partitions for non-transitive pseudo-Anosov flows. It is for this reason that in this note, by adapting Ratner's proof, I prove that Ratner's theorem can be generalized for every expansive flow in dimension 3. Thanks to this result, I am very optimistic that the results of \cite{preprint} can be generalized in the case of non-transitive pseudo-Anosov flows, thus opening a path towards a complete classification of all (transitive and non-transitive) pseudo-Anosov flows in dimension 3. 
\begin{defi}\label{d.expansiveflow}
    Let $M$ be a smooth, closed $3$-manifold endowed with a distance $d$ given by some Riemannian metric. Consider $X^t$ a non-singular $C^0$-flow on $M$. We will say that $X^t$ is \emph{expansive} if there exist $\epsilon, \eta>0$ such that for any $x,y\in M$ and $h:\mathbb{R}\rightarrow \mathbb{R}$ an increasing  homeomorphism with $h(0)=0$ $$ \forall t\in \mathbb{R}~ d(X^t(x),X^{h(t)}(y))<\eta \implies \exists |t_0|<\epsilon ~~ y=X^{t_0}(x)$$

\end{defi}
Classical examples of expansive flows in dimension 3 include transitive or non-transitive Anosov or pseudo-Anosov flows  (see for instance \cite{pseudoanosov}). Conversely, as a result of Theorem 1.5 of \cite{Inaba} and Lemma 2.7 of \cite{Oka}, any expansive flow in dimension 3 is \emph{pseudo-Anosov}: 
\begin{defi}\label{d.prongsingularity}
Consider $\mathcal{F}_h$ the foliation by horizontal lines on $\mathbb{C}$, $\mathcal{D}_2$ the euclidean square $\{z\in\mathbb{C}||\text{Re}(z)|<1 \text{, } |{Im}(z)|<1\}$ and $\pi_p(z)=z^p$ for every $p\in \mathbb{N}^*$. Let   $\mathcal{D}_1= \pi_2(\mathcal{D}_2)$ and $\mathcal{D}_p= \pi_p^{-1}(\mathcal{D}_1)$ for any $p\geq 3$. 

The image of $\mathcal{F}_h$ by $\pi_2$ defines a singular foliation $\mathcal{F}^1_h$. Similarly, for every $p\geq 3$ the pre-image of $\mathcal{F}^1_h$ by $\pi_p$ defines a singular foliation, say $\mathcal{F}^p_h$. We will say that $\mathcal{F}^1_h$ has a \emph{$1$-prong singularity at $0$} and that $\mathcal{F}^p_h$ has a \emph{$p$-prong singularity at $0$}.
\end{defi}

\begin{defi}\label{d.folisingular}
   Consider $M$ a closed $3$-manifold. We will say that $\mathcal{F}$ is \emph{a singular codimension one foliation of $M$} if it is a decomposition of $M$ into \emph{regular} and \emph{singular leaves} satisfying the following properties: 
    \begin{enumerate}
    \item there exists a finite number (possibly zero) of singular leaves $L_1,...,L_n$ 
    \item for every $i\in \llbracket 1,n\rrbracket$ there exists  $C_i\subset L_i$ an embedded circle in $M$ such that $$\text{regular leaves}\cup \{L_1-C_1\}\cup...\cup \{L_n-C_n\}$$ defines on $M-\overset{n}{\underset{i=1}{\cup}}C_i$ a $C^0$ codimension one  foliation. 
        \item for every $x\in C_i$, using the notations of Definition \ref{d.prongsingularity}, there exists $U_x$ a neighborhood of $x$ in $M$ and $h:U_x\rightarrow \mathcal{D}_p\times [0,1]$ a homeomorphism verifying $h(x)=0$ and $h(\mathcal{F}\cap U_x)=\mathcal{F}_h^p\times [0,1]$, where $p\geq 3$ or $p=1$ 
    \end{enumerate}
    We will call $C_i$ a \emph{circle $p$-prong singularity of $\mathcal{F}$} and we will denote the set of circle prong singularities of $\mathcal{F}$ as $\text{Sing}(\mathcal{F})$. 
\end{defi}

\begin{thdef}[Inaba, Matsumoto, Oka]
  Let $M$ be a closed, smooth $3$-manifold, $d$ a distance on $M$ given by some Riemannian metric and $X^t$ an expansive flow on $M$. The flow $X^t$ is a \emph{topological pseudo-Anosov flow}. In other words: 
  \begin{enumerate}
      \item $X^t$ preserves a pair of transverse singular codimension one foliations $F^s$, $F^u$ with a finite number of circle $p$-prong singularities with $p\geq 3$
      \item for every $x,y$ in the same leaf of $F^s$ (resp. $F^u$) there exists an increasing homeomorphism $h:\mathbb{R}\rightarrow \mathbb{R}$ such that $$d(X^t(x),X^{h(t)}(y))\underset{t\rightarrow +\infty}{\longrightarrow}0 \text{ (resp. }d(X^t(x),X^{h(t)}(y))\underset{t\rightarrow -\infty}{\longrightarrow}0\text{)}$$
  \end{enumerate}
  We will call the singular foliations $F^s$ and $F^u$, the \emph{(weak) stable} and \emph{(weak) unstable} foliations of $X^t$ respectively.   
\end{thdef}
\subsection*{Some preliminaries on expansive flows}
Consider $\Phi$ an expansive flow on $M^3$. Denote by  $F^s$ and $F^u$ its weak stable and unstable (singular) foliations and by $NW(\Phi) $ its non-wandering set. 

Our goal in this short section is to define the objects that we are going to use in our construction of Markov partitions and also to prove the existence of a finite set of periodic orbits, whose weak stable leaves form a dense set in $M$ (see Proposition \ref{p.densityfinite}). Our proof of the previous proposition relies on several folkloric results in the theory of expansive flows for which we will provide a proof. 

\begin{figure}
 \begin{minipage}[ht]{0.4\textwidth}
    \centering 
     \vspace{0.8cm}
    \hspace{-0.5cm}
    \includegraphics[width=0.28\textwidth]{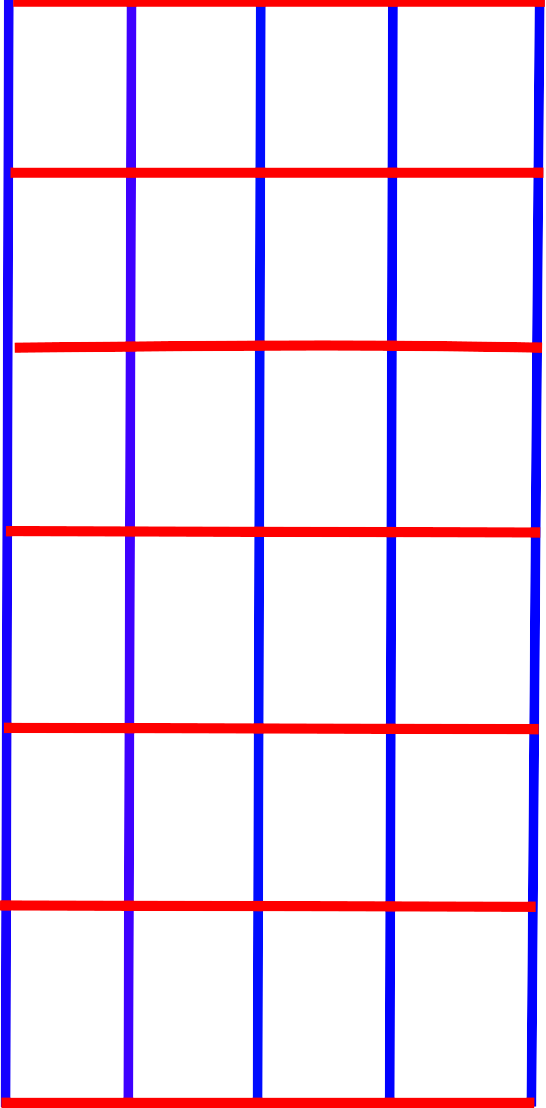}
  \hspace{-1cm}
    \caption*{\quad (a)}
    
  \end{minipage}
 \begin{minipage}[ht]{0.4\textwidth}
 \centering
 \vspace{0.9cm}
    \includegraphics[width=0.55\textwidth]{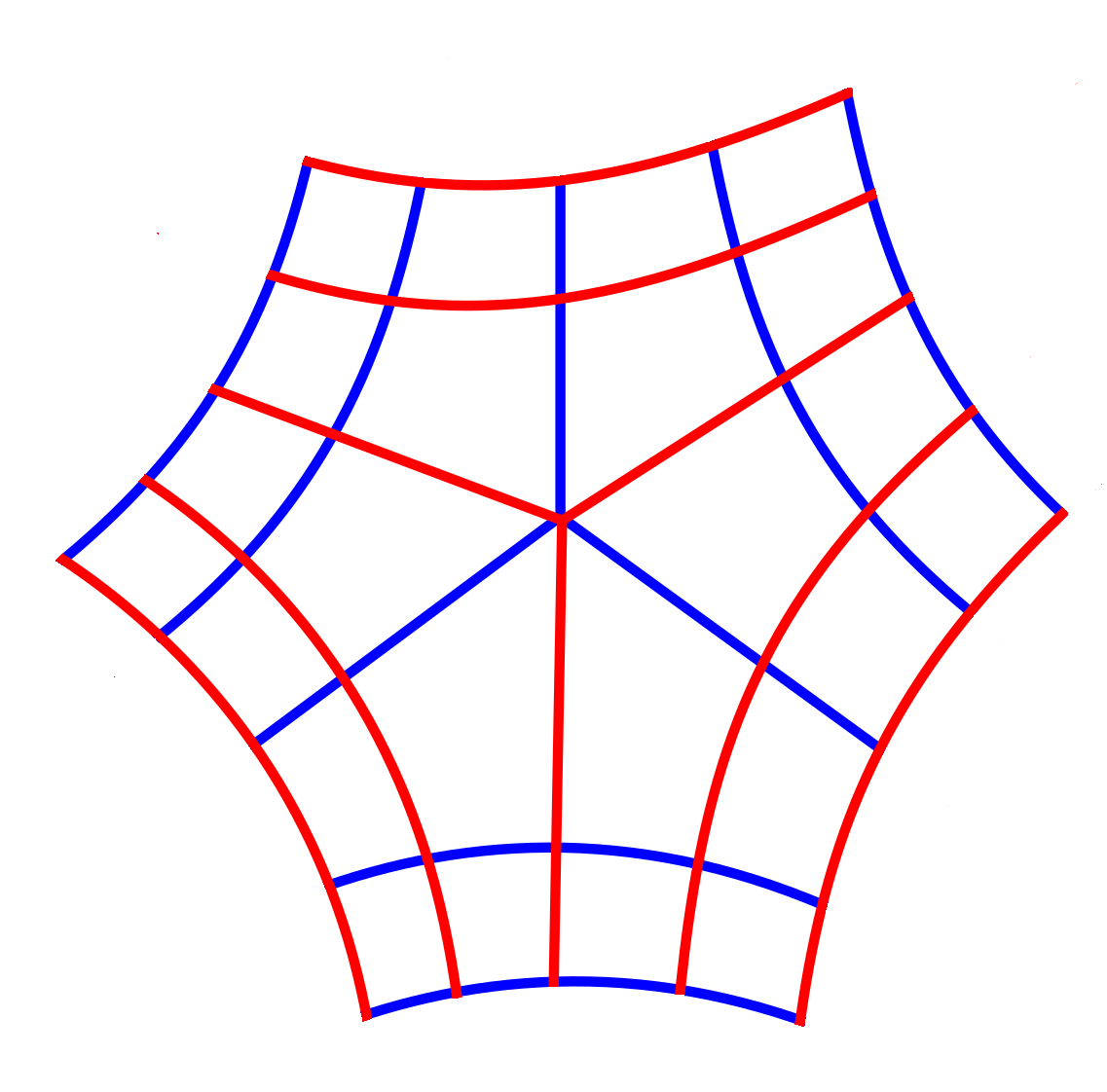}
    
    \caption*{(b)}
    
  \end{minipage}
\caption{(a) A rectangle (b) A standard hexagon}
\label{f.polygons}
  \end{figure}

For every $x\in M$, we are going to denote by $F^s(x)$ (resp. $F^u(x)$) the stable (resp. unstable) leaf containing $x$ and by $F^s_{\epsilon}(x)$ ($\epsilon>0$) the local stable set defined as $$\{y\in F^s(x)| \exists h\in \text{Homeo}_+(\mathbb{R}) ~ h(0)=0 \text{ and } \forall t\geq0 ~ d(\Phi^t(x),\Phi^{h(t)}(y))\leq \epsilon\} $$ We similarly define $F^u_{\epsilon}(x)$. By Theorem 1.5 of \cite{Inaba} and Lemma 2.7 of \cite{Oka}, we have that for every $x\in M$ and every $\epsilon>0$ sufficiently small, $F^s_{\epsilon}(x)$ (resp. $F^u_{\epsilon}(x)$) defines a neighbordhood of $x$ inside $F^s(x)$ (reps. $F^u(x)$) endowed with the leaf topology. 
\begin{defi*}
 A \emph{transverse standard polygon} $P$ in $M$ is an embedded closed topological disk such that: 
 \begin{enumerate}
     \item $P$ is \emph{transverse} to $\Phi$, i.e there exists $\epsilon>0$ such that for every $x\in P$ the segment $\underset{t\in(-\epsilon, \epsilon)}{\cup}\Phi^{t}(x)$ intersects $P$ once
     \item $P$ contains at most one point belonging to a circle prong singularity of $\Phi$
     \item when $\int{P}:=P-\partial P$ intersects a circle $p$-prong singularity, $P$ is bounded by $p$ stable and $p$ unstable segments that alternate between them (see Figure \ref{f.polygons})
     \item when $\int{P}$ does not intersect a circle prong singularity, $P$ is bounded by $2$ stable and $2$ unstable segments that alternate between them. In this case, we will call $P$ a \emph{transverse rectangle}.

\end{enumerate} 
Moreover, for every $x\in P$, the connected component of $F^s(x)\cap P$ (resp. $F^u(x)\cap P$) containing $x$ will be denoted by $P^s(x)$ (resp. $P^u(x)$) and will be called the  \emph{stable} (resp. \emph{unstable}) \emph{segment of $P$ containing $x$}. 
\end{defi*}
\begin{defi*}
    For any transverse rectangle $R$, we will denote by $\partial^{s}R$, $\partial^{u}R$ its stable and unstable boundary and by $\int{R}:=R-(\partial^u R\cup \partial^sR)$ its \emph{interior}. Moreover, any subrectangle $R'$ of $R$ such that $\partial^{s}R'\subset \partial^{s}R$ (resp. $\partial^{u}R'\subset \partial^{u}R$), will be called a \emph{vertical} (resp. \emph{horizontal}) \emph{subrectangle} of $R$. 

\end{defi*}
By Theorem 17A of \cite{Whitney} and the fact that $\Phi$ preserves a pair of transverse singular foliations, it is easy to see that any $x\in M$ is contained in the interior of a transverse standard polygon. In fact, using the compactness of $M$, for any $\epsilon>0$, it is possible to find a finite number of transverse standard polygons $R_1,S_1,R_2,...,R_n,S_n$ of diameter at most $\epsilon$ such that the $S_i$ are pairwise disjoint,  $R_i\subset \int{S_i}$, $R_i$ and $S_i$ have the same number of stable (or unstable) boundary components for every $i$, every orbit segment of length $\alpha>0$ intersects $\underset{i=1}{\overset{n}{\cup}}\int{R_i}$ and every orbit segment of length $4\alpha>0$ intersects at most once every $S_i$ (see also Lemma 2.4 of \cite{Keynes}).

Denote by $p:\underset{i=1}{\overset{n}{\cup}}R_i \rightarrow \underset{i=1}{\overset{n}{\cup}}R_i $ the first return map on $\underset{i=1}{\overset{n}{\cup}}R_i$ and by $B(x,\delta)$ the $\delta$-neighborhood of a point $x\in M$. For a sufficiently small $\rho_0>0$, we have that for any $x\in R_i$ there is a unique continuous map $T:B(x,\rho_0)\rightarrow \mathbb{R}^+$ such that $\Phi^{T(x)}(x)=p(x)\in R_j$ and $\Phi^{T(y)}(y)\in S_j$ for every $y\in B(x,\rho_0)$. We define $p_x(y): B(x,\rho_0) \rightarrow S_j$ as $p_x(y)=\Phi^{T(y)}(y)$. If $p(x)$ and $p_x(y)$ are $\rho_0$-close, then we define $p_x^2(y):=p_{p(x)}(p_x(y))$. We similarly define $p_x^n$ for every $n\in \mathbb{N}^{*}$ and by reversing the orientation of the flow, we can extend our previous definition for all $n\in \mathbb{Z}^{*}$. Notice that for every $\epsilon>0$ sufficiently small and every $x\in R_i$, the function $p_x^n$ is defined on $F^s_{\epsilon}(x)\cap R_i$ for every $n\geq 1$. By Lemma 2.7 of \cite{Oka}, we have that: 


\begin{lemm}\label{l.uniform}
  For any sufficiently small $\epsilon>0$ and every $\delta\in (0,\epsilon)$ there exists $N\in \mathbb{N}$ such that for every $x\in R_i$ and every $m\geq N$ if $p^m(x)\in R_j$, then $p_x^m(F^s_{\epsilon}(x)\cap R_i)\subset F^s_{\delta}(p^m(x))\cap S_j$  
\end{lemm}
Also, as a result of the previous lemma and Lemma 2.5 of \cite{Oka}: 
\begin{lemm}\label{l.foliatedneigh}
    For every sufficiently small $\epsilon>0$, there exists $\eta>0$ such that for any $x,y$ belonging in the same stable (resp. unstable) segment of some $R_i$ and such that $d(x,y)<\eta$, we have that $y\in F^s_{\epsilon}(x)$ (resp. $y\in F^u_{\epsilon}(x)$).
\end{lemm}

\begin{lemm}\label{l.uniformglobal}
    Fix a sufficiently small $\epsilon'>0$. There exists $N\in \mathbb{N}$ such that for every $m\geq N$ (resp. $m\leq -N$) and every $x\in R_i$ there exists a unique continuous function $T: R_i^s(x)\rightarrow \mathbb{R}^{+}$ (resp. $T: R_i^u(x)\rightarrow \mathbb{R}^{-}$) such that $\Phi^{T(x)}(x)=p^m(x)\in R_j$ and $\Phi^{T}(R_i^s(x))\subset F^s_{\epsilon'}(x)\cap S_j$ (resp. $\Phi^{T}(R_i^u(x))\subset F^u_{\epsilon'}(x)\cap S_j$)  .
\end{lemm}
\begin{proof}
    Fix $\epsilon>0$ sufficiently small so that we can apply the Lemmas \ref{l.uniform} and \ref{l.foliatedneigh}. Let $\eta$ be the constant given by Lemma \ref{l.foliatedneigh} associated to $\epsilon$.  Consider a family of $K\in \mathbb{N}$ standard transverse polygons of diameter smaller than $\eta$, whose interiors cover $\underset{i=1}{\overset{n}{\cup}}R_i$. Fix $x\in R_i$. We will prove that by applying positively $\Phi$, the segment $R_i^s$ will be ``uniformly contracted" in the sense of the above statement. The result for $R_i^u$ will follow from a similar argument. 
    
    
     By Lemma \ref{l.foliatedneigh}, there exists a finite number (at most $K$) of points $x_0, x_1,...,x_s\in R_i^s(x)$ such that $R_i^s(x)\subset \underset{i=0}{\overset{n}{\cup}}F^s_{\epsilon}(x_i)$ and $x_{i}\in F^s_{\epsilon}(x_{i-1})$ for every $i\in \llbracket 1,s\rrbracket$. By eventually adding a point to our previous family, we may assume that $x=x_l$, where $l\in \llbracket 0,s\rrbracket$. By Lemma \ref{l.uniform}, if $N$ is sufficiently big there exists a (unique) continuous function $T_N: F^s_{\epsilon}(x)\rightarrow \mathbb{R}^{+}$ such that $\Phi^{T_N(x)}(x)=p^N(x)\in R_j$ and $\Phi^{T_N}(F^s_{\epsilon}(x))\subset F^s_{\epsilon'/K}(p^N(x))\cap S_j$.

    Notice that for every $y\in F^s_{\epsilon}(x)$ the orbit segment $(\Phi^t(y))_{t\in [0,T_N(y)]}$ intersects at least $N$ times $\underset{i=1}{\overset{n}{\cup}}S_i$; hence it intersects at least $r(N)$ times $\underset{i=1}{\overset{n}{\cup}}R_i$, where $r(N)$ depends only on $N$ and $r(N)\underset{N\rightarrow\infty}{\longrightarrow}\infty$. As $x=x_l$ and $x_{l-1},x_{l+1}\in F^s_{\epsilon}(x)$, it follows from Lemma \ref{l.uniform} that for $N$ sufficiently big there exists a (unique) continuous function $T'_N: F^s_{\epsilon}(x_{l-1})\cup F^s_{\epsilon}(x_{l})\cup F^s_{\epsilon}(x_{l+1})\rightarrow \mathbb{R}^{+}$ such that $\Phi^{T'_N(x)}(x)=p^N(x)\in R_j$ and $\Phi^{T'_N(y)}(y)\in F^s_{3\epsilon'/K}(p^N(x))\cap S_j$ for every $y\in F^s_{\epsilon}(x_{l-1})\cup F^s_{\epsilon}(x_{l})\cup F^s_{\epsilon}(x_{l+1})$. We obtain the desired result by induction. 
\end{proof}
\begin{prop}\label{p.weakphase}
   The set $\underset{x\in NW(\Phi)}{\cup}F^s(x)$ is dense in $M$. 
\end{prop}
\begin{proof}
    Assume that this is not the case. As the set $\underset{x\in NW(\Phi)}{\cup}F^s(x)$ contains every singular periodic orbit of $\Phi$, by adding if necessary a rectangle to our family of transverse standard polygons $R_1,S_1,...,R_n,S_n$, we may assume that there exists a transverse rectangle $R_i$ that does not intersect $\underset{x\in NW(\Phi)}{\cup}F^s(x)$. 
    
    Consider $x\in \int{R_i}$ and $y\in \omega(x)\cap \underset{i=1}{\overset{n}{\cup}}R_i$, where $\omega(x)$ denotes the $\omega$-limit of $x$. We will assume for simplicity that $y$ is contained in a transverse rectangle $R_j$ (the general case follows from a similar argument). 
    
    There exists an increasing sequence $t_m\underset{m\rightarrow \infty}{\longrightarrow} +\infty$ such that $x_m:=\Phi^{t_m}(x)\in R_j$ and $x_m\underset{m\rightarrow \infty}{\longrightarrow}y $.  By Lemma \ref{l.uniformglobal}, we have that for $m$ sufficiently big there exists a continuous function $T: R_j^u(x_m)\rightarrow \mathbb{R}^{-}$ such that $\Phi^{T(x_m)}(x_m)=x$ and $\Phi^{T(z)}(z)\in R_i$ for every $z\in R_j^u(x_m)$. In particular, if $[x_m,y]$ is the unique point inside $R_j^u(x_m)\cap R_j^s(y)$, we get that $\Phi^{T([x_m,y])}([x_m,y])\in R_i$ and thus  $F^s(y)$ intersects $R_i$. This leads to an absurd, as by construction  $y\in NW(\Phi)$. 
\end{proof}

\begin{prop}\label{p.densityperiodic}
   The set of periodic orbits of $\Phi$ is dense inside $NW(\Phi)$. 
\end{prop}
\begin{proof}
    Consider $x\in NW(\Phi)$ that does not belong to any periodic orbit of $\Phi$. For every $\epsilon>0$, by 
    changing if necessary our choice of family of transverse standard polygons $R_1,S_1,...,R_n,S_n$, we may assume that $x\in \int{R_i}$, where $R_i$ is a transverse rectangle whose diameter is strictly smaller than $\epsilon$. We would like to show that there exists a periodic orbit of $\Phi$ intersecting $R_i$. 
    
     Since for any $s>0$ the set $\underset{t\in[-s,s]}{\cup}\Phi^t(\int{R_i})$ forms a neighborhood of $x\in NW(\Phi)$, we have that the positive orbit of $R_i$ intersects infinitely many times itself. Hence, there exist $x_m\in R_i$ and $t_m\underset{m\rightarrow \infty}{\longrightarrow} +\infty$ such that $y_m:=\Phi^{t_m}(x_m)\in R_i$. By enlarging very slightly $R_i$ inside $S_i$, we may assume that $x_m,y_m\in \int{R_i}$ and that there exists $\eta>0$ such that $d(x_m,\partial R_i)>\eta$ and $d(y_m,\partial R_i)>\eta$ for every $m$. 
     
     By Lemma \ref{l.uniformglobal}, for $m$ sufficiently big, 
     \begin{itemize}
         \item there exists a (unique) continuous function $T^+_m:R_i^s(x_m)\rightarrow \mathbb{R}^+$ such that $T^+_m(x_m)=t_m$ and $\Phi^{T^+_m}(R_i^s(x_m))\subset F^s_{\eta}(y_m)\cap R_i$
         \item there exists a (unique) continuous function $T^-_m:R_i^u(y_m)\rightarrow \mathbb{R}^-$ such that $T^-_m(y_m)=-t_m$ and $\Phi^{T^-_m}(R_i^u(y_m))\subset F^u_{\eta}(x_m)\cap R_i$
     \end{itemize}
    
    Let $H_{m}$ be the union of all stable segments of $R_i$ that intersect $F^u_{\eta}(x_m)\cap R_i$. The set $H_m$ is a horizontal subrectangle of $R_i$. By our previous arguments, by taking $m$ sufficiently big and by eventually reducing the vertical size of $H_{m}$, we have that there exists a continuous function $T:H_m\rightarrow \mathbb{R}^+$ such that $T(x_m)=t_m$, $\Phi^{T}(H_m)\subset S_i$ and $\Phi^{T}(H_m)\cap R_i$ is a vertical subrectangle of $R_i$. Notice that $\Phi^T$ has a non-zero index along $\partial H_m$. This implies that $\Phi^T$ contains a fixed point inside $H_m$ and thus there exists a periodic orbit of $\Phi$ intersecting $H_m\subset R_i$, which gives us the desired result.  

\end{proof}
\begin{prop}\label{p.densityfinite}
There exists a finite family $\Gamma$ of periodic orbits of $\Phi$ containing every singular periodic orbit of $\Phi$ and such that $\underset{\gamma\in\Gamma}{\cup} F^s(\gamma)$ is dense in $M$. 
\end{prop}

\begin{proof}[Proof of Proposition \ref{p.densityfinite}]
    Take $\mathcal{O}$ a regular periodic orbit of $\Phi$. Consider $U_\mathcal{O}$ the closure of $F^s(\mathcal{O})\cap NW(\Phi)$. It is clear that $U_p$ is a non-empty $\Phi$-invariant closed subset of $NW(\Phi)$. Let us show that it contains a neighborhood of $\mathcal{O}$ inside $NW(\Phi)$. 
    
    Indeed, since $\mathcal{O}$ is a regular periodic orbit and $F^s,F^u$ are transverse, for any $x\in \mathcal{O}$, there exists a rectangle $R_x$ in $M$ transverse to $\Phi$ containing $x$. Fix such a point $x$ and a rectangle $R_x$. For any $z,z'\in R_x$ denote by $[z,z']$ the unique point in $R_x^s(z)\cap R_x^u(z')$. 
    
    Consider $\mathcal{O}'$ a periodic orbit of $\Phi$ intersecting $R_x$ at $y\in \mathcal{O}'$. Recall that $[x,y]\in F^s(\mathcal{O})\cap NW(\Phi)$ (see for instance Lemma 7.2 of \cite{Sm}). The negative orbit of $[x,y]$ accumulates to $\mathcal{O}'$ and since $U_\mathcal{O}$ is $\Phi$-invariant and closed, we have that $\mathcal{O}'\subset U_\mathcal{O}$. This proves that any periodic orbit intersecting $R_x$ belongs in $U_\mathcal{O}$ and thus by Proposition \ref{p.densityperiodic}, we get that $U_\mathcal{O}$ contains a neighborhood of $\mathcal{O}$ inside $NW(\Phi)$. 

    One can prove the same result for a singular periodic orbit by using a transverse standard polygon and the fact that any such polygon is the union of a finite number of transverse rectangles. It follows that since $NW(\Phi)$ is compact, one can find a finite family $\Gamma$ of periodic orbits of $\Phi$, containing all singular periodic orbits of $\Phi$ and such that $\underset{\gamma\in\Gamma}{\cup} \clos{F^s(\gamma)\cap NW(\Phi)}=NW(\Phi)$. We get the desired result by Proposition \ref{p.weakphase}.
\end{proof}

\subsection*{On the construction of Markov partitions}

\begin{defi*}\label{d.markovpartition}
A \emph{Markov partition} of $\Phi$ is a finite family of pairwise disjoint transverse rectangles $R_1,...,R_m$ in $M$ such that: 
\begin{enumerate}
\item Pushing positively by the flow, the first return on $\underset{i=1}{\overset{m}{\cup}}R_i$ of any point $x\in  \underset{i=1}{\overset{m}{\cup}}R_i$ is well defined and will be denoted by $f(x)$. Furthermore, there exists $T>0$ such that for all $x\in M$ there exists $t\in [0,T]$ for which $\Phi^t(x) \in \underset{i=1}{\overset{m}{\cup}}R_i$. 
\item For any two $i,j$ the closure of each connected component of  $f(\int{R_i})\cap \int{R_j}$  (the previous set can be empty) is a vertical subrectangle of $R_j$. 
\item For any two $i,j$ the closure of each connected component of  $f^{-1}(\int{R_i})\cap \int{R_j}$ (the previous set can be empty) is a horizontal subrectangle of $R_j$. 
\end{enumerate}
\end{defi*}

We will use the following classical result in order to construct Markov partitions for $\Phi$: 
\begin{lemm}\label{l.markov}
    Let $R_1,...,R_m$ be a collection of pairwise disjoint rectangles in $M$ that are transverse to $\Phi$ and that intersect every positive and negative orbit of $\Phi$ in a uniformly bounded time. Denote by $f$ the first return map on $\underset{i=1}{\overset{m}{\cup}}R_i$. The family $R_1,...,R_m$ is a Markov partition of $\Phi$ if and only if $$ f(\underset{i=1}{\overset{m}{\cup}}\partial^sR_i)\cap \underset{i=1}{\overset{m}{\cup}}\int{R_i}=\emptyset \text{ and } f^{-1}(\underset{i=1}{\overset{m}{\cup}}\partial^uR_i)\cap \underset{i=1}{\overset{m}{\cup}}\int{R_i}=\emptyset$$
\end{lemm}

\begin{thm*} \label{t.existenceofmarkovpartitions}
Let $\Phi$ be a non-transitive pseudo-Anosov flow on a closed 3-manifold $M$ and $F^s,F^u$ its weak stable and unstable foliations. For any finite family $\Gamma$ of periodic orbits of $\Phi$ containing every singular periodic orbit of $\Phi$ and such that $\underset{\gamma\in \Gamma}{\cup} F^s(\gamma)$ and $\underset{\gamma\in \Gamma}{\cup} F^u(\gamma)$ are dense in $M$, there exists a Markov partition of $\Phi$ formed by rectangles, whose stable and unstable boundaries are contained respectively in $\underset{\gamma \in \Gamma}{\cup}F^s(\gamma)$ and $\underset{\gamma \in \Gamma}{\cup}F^u(\gamma)$.
\end{thm*}

We will prove the above statement in two steps:
\begin{itemize}
    \item \textbf{Step 1}: we will construct a family of pairwise disjoint rectangles transverse to $\Phi$, whose stable and unstable boundaries are contained respectively in $\underset{\gamma \in \Gamma}{\cup}F^s(\gamma)$ and $\underset{\gamma \in \Gamma}{\cup}F^u(\gamma)$ and that intersect every positive and negative orbit of $\Phi$ in a uniformly bounded time
    \vspace{0.2cm}
    \item \textbf{Step 2}: after cutting a finite number of times the above family of rectangles we will show that it gives rise to a Markov partition of $\Phi$ \end{itemize}. 
\begin{proof}[Proof of Step 1]
As before, for every $x\in M$ there exists a transverse standard polygon $D_x$ containing $x$ in its interior and such that for some small $\epsilon_x>0$ the flow $\Phi$ inside $B_x:=\underset{t\in(-\epsilon_x, \epsilon_x)}{\cup}\Phi^{t}D_x$ is $C^0$ conjugated to a constant speed vertical flow. Since $\underset{\gamma\in \Gamma}{\cup} F^s(\gamma)$ and $\underset{\gamma\in \Gamma}{\cup} F^u(\gamma)$ are dense, by eventually reducing the size of the $D_x$, we can ask that the stable boundary of the $D_x$ be contained in $\underset{\gamma \in \Gamma}{\cup}F^s(\gamma)$ and the unstable in $\underset{\gamma \in \Gamma}{\cup}F^u(\gamma)$. The family $(\int{B_x})_{x\in M}$ defines an open cover of $M$; we can thus extract from the previous family a finite cover of $M$, say $(B_{x_i})_{i\in I}$.

The family $(D_{x_i})_{i\in I}$ consists of a finite number of transverse standard polygons, whose union intersects every positive or negative orbit of $\Phi$ in a uniformly bounded time. By cutting every $D_{x_i}$ that is not a rectangle along the invariant manifolds of the unique singular periodic orbit intersecting $D_{x_i}$, we can write $D_{x_i}$ as a union of transverse rectangles. We  thus obtain a finite family of transverse rectangles $(D_i)_{i\in \llbracket 1, s \rrbracket}$, whose boundaries are contained in $\underset{\gamma\in \Gamma}{\cup} (F^s(\gamma)\cup F^u(\gamma))$ and whose union intersects every positive or negative orbit of $\Phi$ in a uniformly bounded time. 

Notice that the rectangles $D_i$ are not necessarily pairwise disjoint. Let us now construct a family of rectangles $(R_i)_{i\in J}$ with the desired properties. Define $R_1:=D_{1}$. Next, if $D_{2}\cap R_1= \emptyset$, we set $R_2:=D_{2}$. If not, take $\epsilon>0$ sufficiently small so that the flow inside $\underset{t\in(-\epsilon, \epsilon)}{\cup}\Phi^{t}D_2$ is conjugated to a constant speed vertical flow. Since $R_1$ is transverse to $\Phi$, for every $y\in D_{2}$ the segment $\underset{t\in(-\epsilon, \epsilon)}{\cup}\Phi^{t}y$ intersects $R_1$ a finite number of times. Therefore, using the density of $\underset{\gamma \in \Gamma}{\cup}F^s(\gamma)$ and $\underset{\gamma \in \Gamma}{\cup}F^u(\gamma)$, by cutting $D_{2}$ into smaller rectangles, we get that there exists a family of subrectangles of $D_{2}$, say $D_{2}^1,...,D_{2}^k$, and $t_1,...,t_k\in (-\epsilon,\epsilon)$ such that:
 \begin{itemize}
     \item $D_{2}= \overset{k}{\underset{i=1}{\cup}}D_{2}^i$
     \item the $D_{2}^i$ have pairwise disjoint interiors
     \item $\partial^s D_{2}^i \subset \underset{\gamma \in \Gamma}{\cup}F^s(\gamma)$ and $\partial^u D_{2}^i \subset \underset{\gamma \in \Gamma}{\cup}F^u(\gamma)$
      \item $\Phi^{t_i}D_{2}^i \cap R_1= \emptyset$
     \item the $t_i$ are pairwise distinct 
    
\end{itemize}
 We set $R_2:=\Phi^{t_1}D_{2}^1$,..., $R_{k+1}:=\Phi^{t_k}D_{2}^k$. Notice that since the flow inside $\underset{t\in(-\epsilon, \epsilon)}{\cup}\Phi^{t}D_2$ is conjugated to a constant speed vertical flow, the rectangles $R_2,...,R_{k+1}$ are pairwise disjoint and transverse to $\Phi$. By repeating the previous argument a finite number of times, we construct a family of pairwise disjoint rectangles $R_1,...,R_N$ with the desired properties.  
 
 \end{proof}

 \begin{proof}[Proof of Step 2]
 Consider $R_1,...,R_N$ the family of rectangles constructed in Step 1.     
 
 For every $\gamma\in \Gamma$, let $K_\gamma^s,K_\gamma^u$ denote two compact subsets of $F^s(\gamma)$ and $F^u(\gamma)$ respectively such that $K_\gamma^s$ is positively invariant by $\Phi$, $K_\gamma^u$ is negatively invariant by $\Phi$,  $\underset{i\in \llbracket 1, N \rrbracket}{\cup}\partial^s R_i\subset \underset{\gamma\in \Gamma}{\cup}K_\gamma^s$ and $\underset{i\in \llbracket 1, N \rrbracket}{\cup}\partial^u R_i\subset \underset{\gamma\in \Gamma}{\cup}K_\gamma^u$. Notice that since the set $\underset{\gamma\in \Gamma}{\cup}K_\gamma^s$ is compact and tangent to $\Phi$, it intersects every $\int{R_i}$ along a finite number of segments $(s^i_j)_{j\in \llbracket 1, J(i)\rrbracket }$. Consider $(l^i_j)_{j\in \llbracket 1, J(i)\rrbracket}$ the connected components of $\underset{\gamma\in \Gamma}{\cup}F^s(\gamma)\cap R_i$ containing the segments $(s^i_j)_{j\in \llbracket 1, J(i)\rrbracket}$. Cut each $R_i$ along the segments $(l^i_j)_{j\in \llbracket 1, J(i)\rrbracket}$ in order to produce $J(i)+1$ new rectangles $R_{i1},...,R_{i(J(i)+1)}$. Repeat the same procedure with $\underset{\gamma\in \Gamma}{\cup}K_\gamma^u$ and the finite family of rectangles $\{R_{ij}|i\in I, j\in \llbracket 1, J(i)+1 \rrbracket \}$. At the end of this procedure, we obtain a new finite family of rectangles in $M$, say $(S_k)_{k\in K}$, that are transverse to $\Phi$ and whose stable and unstable boundaries are contained in $\underset{\gamma\in \Gamma}{\cup}F^s(\gamma)$ and $\underset{\gamma\in \Gamma}{\cup}F^u(\gamma)$ respectively. 
 
 By pushing a little bit certain $S_k$ along the flow, we may assume that the $S_k$ are pairwise disjoint. Furthermore, since the $S_k$ have been obtained by cutting the rectangles $R_i$ it remains true that for every point $x\in M$, the positive and negative orbits of $x$ by $\Phi$ intersect $\underset{k\in K}{\cup}S_k$ in a uniformly bounded time. Therefore, the first return map $f:\underset{k\in K}{\cup}S_k\rightarrow \underset{k\in K}{\cup}S_k$ is well defined. We will prove that $(S_k)_{k\in K}$ is a Markov partition of $\Phi$ thanks to Lemma \ref{l.markov}. Assume that there exist $i,j\in K$ and $x\in\partial^s S_i$ such that $f(x)\in \int{S_j}$. This would imply that the positive orbit by $\Phi$ of the stable boundary of $S_i$ intersects the interior of $S_j$. However, this is impossible since $\partial^s S_i \subset \underset{\gamma\in \Gamma}{\cup}K_\gamma^s$, $\underset{\gamma\in \Gamma}{\cup}K_\gamma^s$ is positively invariant by $\Phi$ and by our construction $\underset{\gamma\in \Gamma}{\cup}K_\gamma^s$ does not intersect the interior of any rectangle in $(S_k)_{k\in K}$. Therefore, $ f(\underset{i=1}{\overset{m}{\cup}}\partial^sR_i)\cap \underset{i=1}{\overset{m}{\cup}}\int{R_i}=\emptyset \text{ and by a similar argument } f^{-1}(\underset{i=1}{\overset{m}{\cup}}\partial^uR_i)\cap \underset{i=1}{\overset{m}{\cup}}\int{R_i}=\emptyset$.
\end{proof}

\end{document}